\documentclass{amsart}
\usepackage{amscd,amssymb,hyperref}

\theoremstyle{plain}
\newtheorem{prop}[subsection]{Proposition}
\newtheorem{thm}[subsection]{Theorem}

\newtheorem{cor}[subsection]{Corollary}

\theoremstyle{remark}
\newtheorem{rem}[subsection]{Remark}

\numberwithin{equation}{section}
\renewcommand{\labelenumi}{{\bfseries\arabic{enumi}.}}

\newcommand{\A}{{\mathcal A}}
\newcommand{\FF}{{\mathcal F}}
\newcommand{\LL}{{\mathcal L}}
\newcommand{\Z}{{\mathbb Z}}
\newcommand{\C}{{\mathbb C}}
\newcommand{\CP}{{\mathbb{CP}}}
\newcommand{\T}{{({\mathbb C}^*)^n}}
\newcommand{\bS}{{\mathbb S}}
\newcommand{\F}{{\mathbb F}}
\newcommand{\bl}{{\boldsymbol{\lambda}}}
\newcommand{\sfB}{{\sf B}}
\newcommand{\sfb}{{\sf b}}
\newcommand{\sfM}{{\sf M}}
\newcommand{\D}{{\Delta}}
\newcommand{\la}{{\lambda }}
\newcommand{\bul}{{\bullet }}

\renewcommand{\b}[1]{\mathbf{#1}}
\renewcommand{\a}{{\alpha }}
\renewcommand{\c}{{\gamma }}
\renewcommand{\L}{{\Lambda }}
\renewcommand{\ll}{{\ell }}

\DeclareMathOperator{\rank}{rank}
\DeclareMathOperator{\codim}{codim}
\DeclareMathOperator{\ii}{i}
\DeclareMathOperator{\id}{id}
\DeclareMathOperator{\Mat}{Mat}
\DeclareMathOperator{\Aut}{Aut}
\DeclareMathOperator{\End}{End}
\DeclareMathOperator{\Hom}{Hom}
\DeclareMathOperator{\TT}{T}
\DeclareMathOperator{\GL}{GL}

\begin{document}

\title[Gauss-Manin Connections]
{Gauss-Manin Connections for Arrangements}
\author[D.~Cohen]{Daniel C.~Cohen$^\dag$}
\address{Department of Mathematics, Louisiana State University,
Baton Rouge, LA 70803}
\email{\href{mailto:cohen@math.lsu.edu}{cohen@math.lsu.edu}}
\urladdr{\href{http://www.math.lsu.edu/~cohen/}
{www.math.lsu.edu/\~{}cohen}}
\thanks{{$^\dag$}Partially supported by
Louisiana Board of Regents grants LEQSF(1996-99)-RD-A-04
and LEQSF(1999-2002)-RD-A-01, and
by National Security Agency grant MDA904-00-1-0038}

\author[P.~Orlik]{Peter Orlik}
\address{Department of Mathematics, University of Wisconsin,
Madison, WI 53706}
\email{\href{mailto:orlik@math.wisc.edu}{orlik@math.wisc.edu}}

\subjclass[2000]{32S22, 14D05, 52C35, 55N25}
% 32S22 Relations with arrangements of hyperplanes
%       (Several complex variables and analytic spaces; Singularities)
% 14D05 Structure of families (Picard-Lefschetz,  monodromy, etc.)
%       (Algebraic geometry; Families, fibrations)
% 52C35 Arrangements of points, flats, hyperplanes
%       (Convex and discrete geometry; Discrete geometry)
% 55N25 Homology with local coefficients, equivariant cohomology
%       (Algebraic topology; Homology and cohomology theories)

\keywords{hyperplane arrangement,
local system, Gauss-Manin connection}

\begin{abstract}
We construct a formal connection on the Aomoto complex of an
arrangement of hyperplanes, and use it to study the Gauss-Manin
connection for the moduli space of the arrangement in the
cohomology of a complex rank one local system.  We prove that the
eigenvalues of the Gauss-Manin connection are integral linear
combinations of the weights which define the local system.
\end{abstract}

\maketitle

\section{Introduction}
\label{sec:intro}

Let $\A=\{H_1,\dots,H_n\}$ be a hyperplane arrangement in $\C^\ll$,
with complement $M=M(\A)=\C^\ll\setminus\bigcup_{j=1}^n H_j$.  Let
$\bl=(\la_1,\dots,\la_n)\in\C^n$ be a collection of weights.
Associated to $\bl$, we have a rank one representation
$\rho:\pi_1(M)\to\C^*$ given by $\c_j\mapsto t_j=\exp(-2\pi\ii\la_j)$
for any meridian loop $\c_j$ about the hyperplane $H_j$ of $\A$, and
an associated rank one local system $\LL$ on $M$.  The need to
calculate the local system cohomology $H^*(M;\LL)$ arises in various
contexts.  For instance, such local systems may be used to study the
Milnor fiber of the non-isolated hypersurface singularity at the
origin obtained by coning the arrangement, see \cite{CS1, CO2}.  In
mathematical physics, local systems on complements of arrangements
arise in the Aomoto-Gelfand theory of multivariable hypergeometric
integrals \cite{AK,Gel1,OT2} and the representation theory of Lie
algebras and quantum groups.  These considerations lead to solutions
of the Knizhnik-Zamolodchikov differential equation from conformal
field theory, see \cite{SV,Va}.  Here, a central problem is the
determination of the Gauss-Manin connection on $H^*(M(\A);\LL)$ for
certain arrangements, and local systems arising from certain weights.
In this paper, we study the Gauss-Manin connection on $H^*(M(\A);\LL)$
for all arrangements, and local systems arising from arbitrary
weights.

The arrangements which arise in the context of the K-Z equations are
the discriminantal arrangements of Schechtman-Varchenko \cite{SV}.
For these arrangements, and certain weights, the Gauss-Manin
connection on $H^{*}(M(\A);\LL)$ was determined by Aomoto \cite{Ao}
and Kaneko \cite{JK}.  The monodromy corresponding to this connection
is a representation of the fundamental group of the moduli space of
all these arrangements, which is a (classical) configuration space,
see \cite{Va}.  Moduli spaces of arbitrary arrangements with a fixed
combinatorial type were defined and investigated by Terao \cite{T1}.
He identified the moduli spaces of certain arrangements, and
determined the Gauss-Manin connection for certain weights.  A priori,
the eigenvalues of this connection are rational functions of the
weights.  Terao found that the eigenvalues are, in fact, integral
linear combinations of the weights.  He asked if this is always the
case.  In Theorem \ref{thm:main}, we prove that the eigenvalues of the
Gauss-Manin connection are indeed integral linear combinations of the
weights for all arrangements and all weights.

Fix the combinatorial type of an arrangement $\A$, and let $\sfB$ be a
smooth, connected component of the moduli space of arrangements of
type $\A$.  There is a fiber bundle $\pi:\sfM \to \sfB$ over $\sfB$.
The fibers of this bundle, $\pi^{-1}(\sfb)=M(\A_\sfb)$, are
complements of arrangements $\A_\sfb$ combinatorially equivalent to
$\A$, so are diffeomorphic to $M(\A)$ (since $\sfB$ is connected).
Given weights $\bl$, we used stratified Morse theory in \cite{CO1} to
construct a complex which computes $H^{*}(M(\A);\LL)$.  In fact, we
constructed a universal complex $(K^{\bul}_{\L}(\A),\D^{\bul}(\b{x}))$
with the property that the specialization $x_j \mapsto
t_j=\exp(-2\pi\ii\la_j)$ calculates $H^*(M(\A);\LL)$.  Here,
$\b{x}=(x_1,\ldots,x_n)$ are the coordinate functions on the complex
$n$-torus $\T$, and $\L=\C[x_1^{\pm 1},\dots,x_n^{\pm 1}]$ is the
coordinate ring.  This construction is reviewed in
Section~\ref{sec:complexes}.

At $\sfb\in \sfB$, we have the corresponding universal complex
$(K^{\bul}_{\L}(\A_{\sfb}),\D^{\bul}(\b{x}))$, its specialization
$(K^{\bul}(\A_{\sfb}),\D^{\bul}(\b{t}))$ and the cohomology of the
latter.  Loops in $\sfB$ based at $\sfb$ induce automorphisms of all
these objects, and hence yield representations of $\pi_1(\sfB,\sfb)$.
In particular, there is a {\em universal representation}
$\pi_1(\sfB,\sfb) \to \Aut_\L K^{\bul}_{\L}(\A_{\sfb})$.  Let
$\b{y}=(y_{1},\dots,y_{n})$ be the coordinates of $\TT_{\b{1}}\T =
\C^{n}$, the holomorphic tangent space of $\T$ at the identity element
$\b{1}=(1,\dots,1)\in\T$.  The exponential map $\TT_{\b{1}}\T \to \T$
is induced by $\exp:\C\to\C^{*}$, $y_{j}\mapsto e^{y_{j}} = x_{j}$.
We call the formal logarithm associated to the universal
representation the {\em formal connection}.

Since the complex $(K^{\bul}(\A_{\sfb}),\D^{\bul}(\b{t}))$ computes
the cohomology of the local system $\LL$ on $M(\A_\sfb)$ corresponding
to the weights $\bl$, the representation $\pi_1(\sfB,\sfb) \to \Aut_\C
H^*(M(\A_\sfb);\LL)$ is induced by the representation
$\pi_1(\sfB,\sfb) \to \Aut_\C K^{\bul}(\A_{\sfb})$.  We realize the
latter as the specialization at $\b{t}$ of the universal
representation.  Similarly, the specialization $y_j \mapsto \la_j$ of
formal connection induces the Gauss-Manin connection on the local
system cohomology.  Given a loop $\c \in \pi_1(\sfB,\sfb)$, we show in
Section~\ref{sec:reps} that the eigenvalues of the corresponding
universal representation matrix are monomials in the $x_j$ with
integer exponents.  In Section~\ref{sec:GM}, we show that the
eigenvalues of the corresponding formal connection matrix are linear
forms in the $y_j$ with integer coefficients.  It follows that the
eigenvalues of the corresponding Gauss-Manin connection matrix in
local system cohomology are integral linear combinations of the
weights, answering Terao's question affirmatively for all arrangements
and all~weights.

The formal connection may be viewed as a connection on a
combinatorial object, the Aomoto complex.  We use the notation and
results of \cite{OT1,OT2}.  Let $A=A(\A)$ be the Orlik-Solomon
algebra of $\A$ generated by the 1-dimensional classes $a_j$,
$1\leq j\leq n$.  It is the quotient of the exterior algebra
generated by these classes by a homogeneous ideal, hence is a
finite dimensional graded $\C$-algebra. There is an isomorphism of
graded algebras $H^*(M;\C) \simeq A(\A)$. For weights $\bl$, the
Orlik-Solomon algebra is a cochain complex with differential given
by multiplication by $a_\bl=\sum_{j=1}^{n}\la_j\,a_j$.  The Aomoto
complex $(A^{\bul}_{R}(\A),a_{\b{y}}\wedge)$ is a universal
complex with the property that the specialization $y_j \mapsto
\la_j$ calculates $H^*(A^{\bul},a_\bl\wedge)$.  Here,
$R=\C[y_1,\dots,y_n]$ is the coordinate ring of $\C^n$, the
holomorphic tangent space of $\T$ at $\b{1}$.  In
\cite[Thm.~2.13]{CO1}, we showed that the Aomoto complex
$(A^{\bul}_{R}(\A),a_{\b{y}}\wedge)$ is chain equivalent to the
linearization of the universal complex
$(K^{\bul}_{\L}(\A),\D^{\bul}(\b{x}))$.

Call a system of weights $\bl$ or the corresponding local system
$\LL$ {\em combinatorial} if local system cohomology is
quasi-isomorphic to Orlik-Solomon algebra cohomology,
\[
H^*(M(\A);\LL) \simeq H^*(A^\bul(\A),a_\bl\wedge).
\]
The set of combinatorial weights is open and dense in $\C^n$.  See
\cite{ESV, STV} for sufficient conditions.  For combinatorial weights,
the Gauss-Manin connection in local system cohomology coincides with
that in the cohomology of the Orlik-Solomon algebra.  Thus, if the
eigenvalues of the former are integer linear combinations of the
weights, then so are those of the latter.  In Section \ref{sec:GM2},
we show that, in fact, the eigenvalues of the
combinatorial
Gauss-Manin connection
in Orlik-Solomon algebra cohomology are integer linear combinations of
the weights for all weights.  Since the Aomoto complex
$A^{\bul}_{R}(\A)$ is the linearization of the universal complex
$K^{\bul}_{\L}(\A)$, results on the universal representation on
$K^{\bul}_{\L}(\A)$ inform on the formal connection on
$A^{\bul}_{R}(\A)$, and its specializations, for arbitrary weights.

Call a system of weights $\bl$ or the corresponding local system
$\LL$ {\em non-resonant} if the Betti numbers of $M$ with
coefficients in  $\LL$ are minimal. The set of non-resonant
weights is open and dense in $\C^n$, but does not coincide with
the set of combinatorial weights.  The cohomology of non-resonant
local systems is known.  A detailed account, including sufficient
conditions, is found in \cite{OT2}. For non-resonant weights we
have
\begin{equation} \label{eq:NonResCohomology}
H^q(M;\LL)=0 \mbox{ for } q\neq \ll, \mbox{ and } \dim
H^{\ll}(M;\LL)=|e(M)|,
\end{equation}
where $e(M)$ is the Euler characteristic, see \cite{ESV,STV,Yuz}.  If
the weights are both combinatorial and non-resonant, the Gauss-Manin
connection may be studied effectively by combinatorial means.  In
particular, explicit bases for the single non-vanishing cohomology
group are known, see \cite{FT}.

Several authors have studied Gauss-Manin connections using such a
basis. See Aomoto \cite{Ao} and Kaneko \cite{JK} for
discriminantal arrangements, and Kanarek \cite{HK} for the
connection arising when a single hyperplane in the arrangement is
allowed to move. For general position arrangements, the
Gauss-Manin connection matrices were computed by Aomoto-Kita
\cite{AK}.  These connection matrices were obtained by Terao
\cite{T1} for a larger class of arrangements.  Like the
eigenvalues, the entries of these matrices were known to be
rational functions of the weights.  Both Aomoto-Kita and Terao
found that these entries were, in fact, integer linear
combinations of the weights, and Terao asked if this is the case
in general.  Our work here was motivated by these results and this
question.

\section{Cohomology Complexes}
\label{sec:complexes}

For an arbitrary complex local system $\LL$ on the complement of an
arrangement $\A$, we used stratified Morse theory in \cite{C1} to
construct a complex $(K^\bul(\A),\D^\bul)$, the cohomology of which is
naturally isomorphic to $H^*(M;\LL)$, the cohomology of $M$ with
coefficients in $\LL$.  We now recall this construction in the context
of rank one local systems, and record several related complexes and
relevant results from \cite{C1,CO1}.

Choose coordinates $\b{u}=(u_1,\dots,u_\ll)$ on $\C^\ll$, and let
$\A=\{H_1,\dots,H_n\}$ be a hyperplane arrangement in $\C^\ll$, with
complement $M=M(\A)=\C^\ll\setminus \bigcup_{j=1}^n H_j$.  We assume
throughout that $\A$ contains $\ll$ linearly independent hyperplanes.
For each $j$, let $f_j$ be a linear polynomial with $H_j = \{\b{u} \in
\C^\ll \mid f_j(\b{u})=0\}$.  Let $\bl=(\la_1,\dots,\la_n)\in\C^n$ be
a system of weights.  Associated to $\bl$, we have
\begin{enumerate}
\item[(1)] a flat connection on the trivial line bundle over $M$, with
connection form $\nabla=d+\omega_{\bl}\wedge:\Omega^0 \to \Omega^1$,
where $d$ is the exterior differential operator with respect to the
coordintes $\b{u}$, $\omega_{\bl} = \sum_{j=1}^n \la_j\,d\log(f_j)$,
and $\Omega^q$ is the sheaf of germs of holomorphic differential forms
of degree $q$ on $M$;

\item[(2)] a rank one representation $\rho:\pi_1(M) \to \C^*$, given
by $\rho(\c_j) =t_j$, where $\b{t}=(t_1,\dots,t_n)\in (\C^*)^n$ is
defined by $t_j= \exp(-2\pi\ii\la_j)$, and $\c_j$ is any meridian loop
about the hyperplane $H_j$ of $\A$; and

\item[(3)] a rank one local system $\LL=\LL_{\b{t}}=\LL_{\bl}$ on $M$
associated to the representation $\rho$ (resp., the flat connection
$\nabla$).
\end{enumerate}
Note that weights $\bl$ and $\bl'$ yield identical representations
and local systems if $\bl-\bl'\in\Z^n$.

\begin{rem} \label{rem:strat}
The arrangement $\A$ determines a Whitney stratification of $\C^\ll$,
with codimension zero stratum given by the complement $M$.  To
describe the strata of higher codimension, recall that an edge of $\A$
is a nonempty intersection of hyperplanes.  Associated to each
codimension $p$ edge $X$, there is a stratum $S_X=X \setminus \bigcup
Y$, where the union is over all edges $Y$ of $\A$ which satisfy
$Y\subsetneq X$.  Note that $S_X=M(\A^X)$ may be realized as the
complement of the arrangement $\A^X$ in $X$, see \cite{OT1}.
\end{rem}

Let $\FF$ be a complete flag (of affine subspaces) in $\C^\ll$,
\begin{equation} \label{eq:flag}
\FF:\quad \emptyset = \FF^{-1} \subset \FF^0 \subset \FF^1 \subset
\FF^2 \subset\dots \subset \FF^\ll = \C^\ll,
\end{equation}
transverse to the stratification determined by $\A$, so that $\dim
\FF^q\cap S_X = q-\codim S_X$ for each stratum, where a negative
dimension indicates that $\FF^q\cap S_X=\emptyset$.  For an explicit
construction of such a flag, see \cite[\S1]{C1}.  Let $M^{q} = \FF^q
\cap M$ for each $q$.  Let $K^q=H^q(M^q,M^{q-1};\LL)$, and denote by
$\D^q$ the boundary homomorphism $H^{q}(M^{q},M^{q-1};\LL) \to
H^{q+1}(M^{q+1},M^{q};\LL)$ of the triple $(M^{q+1},M^q,M^{q-1})$.
The following compiles several results from \cite{C1}.

\begin{thm}\label{thm:Kdot}
Let $\LL$ be a complex rank one local system on the complement $M$ of
an arrangement $\A$ in $\C^\ll$.

\begin{enumerate}
\item \label{item:Kdot1}
For each $q$, $0\le q \le \ll$, we have $H^i(M^q,M^{q-1};\LL) = 0$ if
$i \neq q$, and $\dim_\C H^q(M^q,M^{q-1};\LL) = b_q(\A)$ is equal to
the $q$-th Betti number of $M$ with trivial local coefficients $\C$.

\item \label{item:Kdot2}
The system of complex vector spaces and linear maps
$(K^\bul,\D^\bul)$,
\[
K^0 \xrightarrow{\ \D^{0}\ } K^1 \xrightarrow{\ \D^1\ } K^2
\xrightarrow{\phantom{\D^{1}}} \cdots \xrightarrow{\phantom{\D^{1}}}
K^{\ll-1} \xrightarrow{\ \D^{\ll-1}\,} K^\ll,
\]
is a complex $(\D^{q+1}\circ\D^q=0)$.  The cohomology of this complex
is naturally isomorphic to $H^*(M;\LL)$, the cohomology of $M$ with
coefficients in $\LL$.
\end{enumerate}
\end{thm}

The dimensions of the terms, $K^{q}$, of the complex
$(K^\bul,\D^\bul)$ are independent of $\b{t}$ (resp., $\bl$, $\LL$).
Write $\D^\bul=\D^\bul(\b{t})$ to indicate the dependence of the
complex on $\b{t}$, and view these boundary maps as functions of
$\b{t}$.  Let $\L=\C[x_1^{\pm 1},\dots,x_n^{\pm 1}]$ be the ring of
complex Laurent polynomials in $n$ commuting variables, and for each
$q$, let $K^q_\L= \L \otimes_{\C} K^q$.

\begin{thm}[{\cite[Thm.~2.9]{CO1}}] \label{thm:univcx}
For an arrangement $\A$ of $n$ hyperplanes with complement $M$, there
exists a universal complex $(K^\bul_\L,\D^\bul(\b{x}))$ with the
following properties:
\begin{enumerate}
\item \label{item:univcx1}
The terms are free $\L$-modules, whose ranks are given by the Betti
numbers of $M$, $K^q_\L \simeq \L^{b_q(\A)}$.

\item \label{item:univcx2}
The boundary maps, $\D^q(\b{x}): K^q_\L \to K^{q+1}_\L$ are
$\L$-linear.

\item \label{item:univcx3}
For each $\b{t}\in\T$, the specialization $\b{x} \mapsto \b{t}$ yields
the complex $(K^\bul,\D^\bul(\b{t}))$, the cohomology of which is
isomorphic to $H^*(M;\LL_\b{t})$, the cohomology of $M$ with
coefficients in the local system associated to $\b{t}$.
\end{enumerate}
\end{thm}

The entries of the boundary maps $\D^{q}(\b{x})$ are elements of the
Laurent polynomial ring $\L$, the coordinate ring of the complex
algebraic $n$-torus.  Via the specialization $\b{x} \mapsto \b{t} \in
\T$, we view them as holomorphic functions $\T\to\C$.  Similarly, for
each $q$, we view $\D^{q}(\b{x})$ as a holomorphic map
$\D^{q}:\T\to\Mat(\C)$, $\b{t}\mapsto \D^{q}(\b{t})$.

\begin{rem} \label{rem:KdotAtOne}
If $\b{t}=\b{1}$ is the identity element of $\T$, the associated local
system $\LL_\b{1}$ is trivial.  Consequently, the specialization
$\b{x} \mapsto \b{1}$ yields a complex $(K^\bul,\D^\bul(\b{1}))$ whose
cohomology gives $H^*(M;\C)$.  Since $\dim K^q = b_q(\A) = \dim
H^q(M;\C)$, the boundary maps of this complex are necessarily trivial,
$\D^q(\b{1})=0$ for each $q$.
\end{rem}

There is an analogous universal complex for the cohomology,
$H^{*}(A^{\bul},a_{\bl}\wedge)$, of the Orlik-Solomon algebra
$A=A(\A)$.  Let $R=\C[y_{1},\dots,y_{n}]$ be the polynomial ring.  The
{\em Aomoto complex} $(A^{\bul}_{R},a_{\b{y}}\wedge)$ has terms
$A^{q}_{R}=R\otimes_{\C} A^{q}\simeq R^{b_q(\A)}$, and boundary maps
given by $p(\b{y})\otimes\eta \mapsto \sum y_{j}p(\b{y}) \otimes a_{j}
\wedge \eta$.  For $\bl \in \C^{n}$, the specialization $\b{y}\mapsto
\bl$ of the Aomoto complex $(A^{\bul}_{R},a_{\b{y}}\wedge)$ yields the
Orlik-Solomon algebra complex $(A^{\bul},a_{\bl}\wedge)$.

A choice of basis for the Orlik-Solomon algebra of $\A$ yields a basis
for each term $A^{q}_{R}$ of the Aomoto complex.  Let $\mu^{q}(\b{y})$
denote the matrix of $a_{\b{y}}\wedge:A^{q}_{R} \to A^{q+1}_{R}$ with
respect to a fixed basis.  The following results were established in
\cite{CO1}.

\begin{thm} \label{thm:approx}
\

\begin{enumerate}
\item \label{item:approx1}
For each $q$, the entries of $\mu^{q}(\b{y})$ are integral linear
forms in $y_1,\dots,y_n$.

\item \label{item:approx2}
For any arrangement $\A$, the Aomoto complex
$(A^{\bul}_{R},\mu^\bul(\b{y}))$ is chain equivalent to the
linearization of the universal complex
$(K^{\bul}_{\L},\D^{\bul}(\b{x}))$.
\end{enumerate}
\end{thm}

\section{Representations}
\label{sec:reps}

Let $\A$ be an arrangement of $n$ hyperplanes in $\C^\ll$ as
above, and let $\sfB$ be a smooth, connected component of the
moduli space of arrangements with the combinatorial type of $\A$.
This moduli space is a locally closed subspace of $(\CP^\ll)^n$.
We refer to \cite{OT2,T1} for the precise definition of this
moduli space, and to Section \ref{sec:example} for an example.  In
this section, we extend the constructions of the previous section
to produce representations of the fundamental group of $\sfB$
related to the cohomology of the complement of $\A$ with
coefficients in a rank one local system.

Denote the coordinates on $(\CP^\ll)^n$ by
$\b{z}=(\b{z}^{1},\dots,\b{z}^{n})$, where
$\b{z}^{i}=(z^{i}_{0}:\dots:z^{i}_{\ll})$, and recall that the
coordinates on $\C^\ll$ are denoted by $\b{u}=(u_1,\dots,u_\ll)$.
There is a fiber bundle $\pi:\sfM \to \sfB$, see \cite[\S3]{T1}.  The
total space may be described as
\[
\sfM = \{(\b{z},\b{u}) \in (\CP^\ll)^n \times \C^\ll \mid \b{z} \in \sfB
\text{ and } \b{u} \in \pi^{-1}(\b{z})\},
\]
and the projection is given by $\pi(\b{z},\b{u})=\b{z}$.  For
$\sfb\in\sfB$, the fiber $\sfM_\sfb = \pi^{-1}(\sfb)$ is the
complement, $\sfM_\sfb=M(\A_\sfb)$, of the arrangement $\A_\sfb$
combinatorially equivalent to $\A$.  The closure,
$\overline{\sfM}_\sfb$, of the fiber is homeomorphic to $\C^\ll$, and
admits a Whitney stratification determined by the arrangement
$\A_\sfb$ as in Remark~\ref{rem:strat}.  Let $\FF_\sfb$ be a flag in
$\overline{\sfM}_\sfb$ that is transverse to $\A_\sfb$ as in
\eqref{eq:flag}.  Evidently, these flags may be chosen to vary
smoothly with $\sfb$.

Recall that the hyperplanes of $\A$ are defined by linear polynomials
$f_j = f_j(\b{u})$.  Since $\sfB$ is by assumption connected, for
every $\sfb \in \sfB$, the arrangement $\A_\sfb$ is lattice-isotopic
to $\A$ in the sense of Randell \cite{Ra}.  Consequently, there are
smooth functions $f_j(\b{z},\b{u})$ on $\sfM$ so that, for each $\sfb
\in \sfB$, the hyperplanes of $\A_\sfb$ are defined by
$f_j(\sfb,\b{u})$.

Given $\b{t}\in\T$ (or weights $\bl\in\C^n$) and $\sfb\in\sfB$, denote
the corresponding local system on $\sfM_\sfb$ by $\LL(\sfb)$.  In this
context, the construction of the previous section yields vector
bundles $\b{K}^q$ over $\sfB$ for $0\le q\le\ll$ as follows.  For
$\sfb\in\sfB$, let $\sfM_{\sfb}^{q}=\FF^q_\sfb\cap\sfM_\sfb^{}$ and
\[
K^q(\sfb)= H^q(\sfM_\sfb^{q},\sfM_\sfb^{q-1}; \LL(\sfb)).
\]
Since $\pi:\sfM\to\sfB$ is locally trivial, the natural projection
$\pi^{q}:\b{K}^{q}\to\sfB$, where $\b{K}^{q}=\bigcup_{\sfb\in\sfB}
K^{q}(\sfb)$, is a vector bundle.  The transition functions of this
vector bundle are locally constant.

If $\c:I\to\sfB$ is a path, then the induced bundle
$\c^{*}(\b{K}^{q})$ is trivial.  Consequently there is a canonical
linear isomorphism $K^{q}(\c(0)) \to K^{q}(\c(1))$, from the fiber
over the initial point of $\c$ to that over the terminal point, which
depends only on the homotopy class of the path.  Fix a basepoint
$\sfb\in\sfB$, and write $K^{\bul}=K^{\bul}(\sfb)$.  The operation of
parallel translation of fibers over curves in $\sfB$ in the vector
bundle $\pi^{q}:\b{K}^{q}\to\sfB$ provides a complex representation of
rank $b_q(\A)$,
\begin{equation} \label{eq:Krep}
\Phi^{q}:\pi_{1}(\sfB,\sfb) \longrightarrow \Aut_\C(K^{q}).
\end{equation}
To indicate the dependence of the representation $\Phi^q$ on
$\b{t}\in\T$, write $\Phi^q=\Phi^q(\b{t})$.

\begin{thm} \label{thm:TrivialAtOne}
If $\b{t}=\b{1}$ is the identity element of $\T$, then the
corresponding representation $\Phi^q(\b{1})$ is trivial for each $q$.
That is, for every $\c\in\pi_1(\sfB,\sfb)$, we have
$\Phi^q(\b{1})(\c)=\id:K^q \to K^q$ for each $q$, $0 \le q \le \ll$.
\end{thm}
\begin{proof}
For the trivial local system $\LL(\sfb)=\C$ associated to
$\b{t}=\b{1}$, the long exact cohomology sequence of the pair
$(\sfM_\sfb^q, \sfM_\sfb^{q-1})$ splits into short exact sequences
\[
0 \to H^i(\sfM_\sfb^q,\sfM_\sfb^{q-1};\C) \to
H^i(\sfM_\sfb^q;\C) \to
H^i(\sfM_\sfb^{q-1};\C) \to 0,
\]
see \cite[III.3]{GM} and \cite[Rem.~5.4]{C1}.  In particular, the
$q$-th relative cohomology group $K^q$ is canonically isomorphic to
$H^q(\sfM_\sfb^q;\C)=H^q(\sfM_\sfb;\C)$, the $q$-th cohomology of
$\sfM_\sfb$ with constant coefficients $\C$, see
Remark~\ref{rem:KdotAtOne}.  So it suffices to show that the
fundamental group of $\sfB$ acts trivially on $H^q(\sfM_\sfb;\C)$.

Let $i_{\sfb}:\sfM_{\sfb} \to \sfM$ denote the inclusion of the fiber
in the total space of the bundle $\pi:\sfM \to \sfB$.  It is known
\cite{Serre} that the image, $i_{\sfb}^{*}H^{*}(\sfM;\C) \subseteq
H^{*}(\sfM_{\sfb};\C)$, of the cohomology of the total space $\sfM$
with (trivial) coefficients in the field $\C$ is invariant under the
action of $\pi_{1}(\sfB,\sfb)$.  Consider the logarithmic forms on
$\sfM$ defined by
\[
\omega_j(\b{z},\b{u}) =
\frac{d_{\b{z}}f_j(\b{z},\b{u}) + d_{\b{u}}f_j(\b{z},\b{u})}
{f_j(\b{z},\b{u})}, \quad \text{where} \quad
d_{\b{y}}g = \sum_{i=1}^k \frac{\partial g}{\partial y_i}\, dy_i
\]
denotes the gradient of $g$ with respect to the variables
$\b{y}=(y_1,\dots,y_k)$.  Clearly these forms represent non-trivial
classes in $H^*(\sfM;\C)$.  Furthermore, we have
\[
i^*_\sfb\omega_j(\b{z},\b{u}) = \omega_j(\sfb,\b{u}) =
\frac{d_{\b{z}}f_j(\sfb,\b{u}) + d_{\b{u}}f_j(\sfb,\b{u})}
{f_j(\sfb,\b{u})}=
\frac{d_{\b{u}}f_j(\sfb,\b{u})}
{f_j(\sfb,\b{u})}=
d_{\b{u}}\log(f_j(\sfb,\b{u})).
\]
As is well known, the forms $\omega_j(\sfb,\b{u})$ generate the
cohomology ring of $\sfM_\sfb=M(\A_\sfb)$.  It follows that
$\sfM_{\sfb}$ is totally nonhomologous to zero in $\sfM$ with respect
to $\C$: The inclusion $i_\sfb:\sfM_\sfb \to \sfM$ induces a
surjection $i^*_\sfb:H^*(\sfM;\C) \to H^*(\sfM_\sfb;\C)$ in cohomology
with trivial coefficients $\C$.  Consequently, the fundamental group
$\pi_1(\sfB,\sfb)$ acts trivially on the $H^{q}(\sfM_{\sfb};\C)$ for
each $q$, and the representation $\Phi^{q}(\b{1})$ is trivial.
\end{proof}

Denote the boundary homomorphism of the triple $(\sfM_\sfb^{q+1},
\sfM_\sfb^{q}, \sfM_\sfb^{q-1})$ in cohomology with local coefficients
$\LL(\sfb)$ determined by $\b{t}$ by
$\D^{q}(\b{t})=\D^{q}_{\sfb}(\b{t}):K^{q} \to K^{q+1}$.

\begin{cor} \label{prop:commute}
For each $\b{t}\in\T$ and each $\c\in\pi_1(\sfB,\sfb)$, the
automorphisms $\Phi^q(\b{t})(\c):K^{q}\to K^{q}$, $0\le q\le \ll$,
comprise a chain automorphism $\Phi^\bul(\b{t})(\c)$ of the complex
$(K^\bul,\D^{\bul}(\b{t}))$.
\end{cor}
\begin{proof}
By Theorem \ref{thm:TrivialAtOne}, the result holds at $\b{t}=\b{1}$.
Therefore it holds for $\b{t}$ close to $\b{1}$.  The result follows.
\end{proof}

We abbreviate the above result by writing
$\Phi^\bul(\b{t}):\pi_1(\sfB,\sfb) \to \Aut_\C(K^\bul)$.  For
$\c\in\pi_1(\sfB,\sfb)$, the automorphism
$\Phi^q(\c)=\Phi^q(\b{t})(\c)$ may be viewed as a holomorphic function
of $\b{t}$:
\[
\Phi^q(\c):\T \to \Aut_\C(K^q),\quad \b{t} \mapsto \Phi^q(\b{t})(\c).
\]
Recall that $\L=\C[x_1^{\pm 1},\dots,x_n^{\pm 1}]$.  Let
$(K^\bul_\L,\D^{\bul}(\b{x}))$ be the universal complex of the
arrangement $\A_\sfb$ from Theorem~\ref{thm:univcx}.  By the
continuity of the functions $\Phi^q(\c)$, we have the following
extension of this result.

\begin{thm} \label{thm:univend}
For each $\c\in\pi_1(\sfB,\sfb)$, there is a chain map
$\Phi^\bul(\b{x})(\c): K^\bul_\L \to K^\bul_\L$ so that the
specialization $\b{x} \mapsto \b{t}$ yields the chain automorphism
$\Phi^\bul(\b{t})(\c)$ of the complex $(K^\bul,\D^\bul(\b{t}))$.
This provides a representation $\Phi^\bul(\b{x}):\pi_1(\sfB,\sfb)
\to \End_\L(K^\bul_\L)$ which specializes to the representation
$\Phi^\bul(\b{t}):\pi_1(\sfB,\sfb) \to \Aut_\C(K^\bul)$.
\end{thm}

Call $\Phi^\bul(\b{x}):\pi_1(\sfB,\sfb) \to \End_\L(K^\bul_\L)$ the
{\em universal representation}.

\begin{thm} \label{thm:univaut}
For each $q$ and each $\c\in\pi_{1}(\sfB,\sfb)$, the eigenvalues of
$\Phi^{q}(\b{x})(\c)$ are monomials functions of the form
$r(\b{x})=x_{1}^{m_{1}} \cdots x_{n}^{m_{n}}$, where $m_{j}\in\Z$.
\end{thm}
\begin{proof}
Given $q$ and $\c\in\pi_{1}(\sfB,\sfb)$, let $r(\b{x})$ be an
eigenvalue of $\Phi^q(\b{x})(\c)$.  Then $r(\b{t})$ is an eigenvalue
of $\Phi^q(\b{t})(\c) \in \Aut_\C(K^{q}) \simeq \GL(b_{q}(\A),\C)$ for
every $\b{t}\in\T$.  It follows that the function $r:\T \to \C$,
$\b{t} \mapsto r(\b{t})$ is single-valued and has no poles.  Thus,
$r(\b{x})$ is a Laurent polynomial in $x_1,\dots,x_n$.  Write
$r(\b{x})= x_{1}^{m_{1}} \cdots x_{n}^{m_{n}} \cdot p(\b{x})$, where
$p(\b{x})$ is a polynomial.  Since $\Phi^q(\b{t})(\c)$ is an
automorphism for every $\b{t}\in\T$, we have $p(\b{t}) \neq 0$ for all
$\b{t}$.  Thus, $p(\b{x})=c\in\C^{*}$ is a non-zero constant, and
$r(\b{x}) = c\cdot x_{1}^{m_{1}} \cdots x_{n}^{m_{n}}$ is a unit in
$\L$.  Using Theorem \ref{thm:TrivialAtOne}, we have $c=1$.
\end{proof}

Thus for every $\c\in\pi_1(\sfB,\sfb)$, the maps $\Phi^q(\b{x})(\c)$
are automorphisms, so we write $\Phi^\bul(\b{x}):\pi_1(\sfB,\sfb) \to
\Aut_\L^{}(K^\bul_\L)$.

\begin{cor} \label{cor:Keigen}
For each $\b{t}\in\T$, the eigenvalues of the automorphism
$\Phi^q(\b{t})(\c)$ are evaluations $r(\b{t})$ of the monomial
functions $r(\b{x})$.
\end{cor}

Given $\b{t}\in\T$ with associated local system $\LL(\sfb)$ on
$\sfM_\sfb$, there are also vector bundles $\b{H}^q\to\sfB$ over the
moduli space, defined by
$\b{H}^q=\bigcup_{\sfb\in\sfB}H^q(\sfM_\sfb;\LL(\sfb))$ for each $q$.
As above, parallel translation of the fibers in this bundle over
curves in the base gives rise to a representation
\[
\Psi^{q}=\Psi^{q}(\b{t}):\pi_{1}(\sfB,\sfb) \longrightarrow
\Aut_\C(H^{q}(\sfM_{\sfb};\LL(\sfb))).
\]

By Theorem \ref{thm:Kdot}.\ref{item:Kdot2}, the cohomology of the
complex $(K^{\bul},\D^{\bul}(\b{t}))$ is naturally isomorphic to
$H^{*}(\sfM_{\sfb};\LL(\sfb))$.  Furthermore, parallel translation
induces the representation $\Phi^{\bul}(\b{t}):\pi_{1}(\sfB,\sfb) \to
\Aut_\C(K^{\bul})$ of \eqref{eq:Krep} on this complex.  Thus by
functoriality, we have the following.

\begin{thm} \label{thm:InducedRep}
The cohomology representation $\Psi^{q}(\b{t}) = \Phi^{q}(\b{t})^{*}$
is induced by~$\Phi^{\bul}(\b{t})$.  In other words, for each
$\c\in\pi_1(\sfB,\sfb)$, the automorphism $\Psi^\bul(\b{t})(\c)$ in
cohomology is induced by the automorphism $\Phi^\bul(\b{t})(\c)$ of
the complex $(K^{\bul},\D^{\bul}(\b{t}))$.
\end{thm}

\begin{cor} \label{cor:Heigen}
For each $\b{t}\in\T$, the eigenvalues of the automorphism
$\Psi^q(\b{t})(\c)$ are evaluations of monomial functions.
\end{cor}

\section{Connections}
\label{sec:GM}

The vector bundles $\b{K}^q\to\sfB$ and $\b{H}^q\to\sfB$ over the
moduli space constructed above support Gauss-Manin connections
corresponding to the representations $\Phi^q(\b{t})$ and
$\Psi^q(\b{t})$ of the fundamental group of $\sfB$.  We now study
these connections in light of the results of the previous section.

Over a manifold such as $\sfB$, there is a well known equivalence
between local systems and complex vector bundles equipped with
flat connections, see \cite{De,Ko}.  Let $\b{V}\to\sfB$ be such a
bundle, with connection $\nabla$.  The latter is a $\C$-linear map
$\nabla:\Omega^0(\b{V}) \to \Omega^1(\b{V})$, where
$\Omega^p(\b{V})$ denotes the complex $p$-forms on $\sfB$ with
values in $\b{V}$, which satisfies $\nabla(f\sigma)=\sigma df +
f\nabla(\sigma)$ for a function $f$ and
$\sigma\in\Omega^0(\b{V})$.  The connection extends to a map
$\nabla:\Omega^p(\b{V}) \to \Omega^{p+1}(\b{V})$ for $p\ge 0$, and
is flat if the curvature $\nabla\circ\nabla$ vanishes.  Call two
connections $\nabla$ and $\nabla'$ on $\b{V}$ isomorphic if
$\nabla'$ is obtained from $\nabla$ by a gauge transformation,
$\nabla'=g\circ\nabla\circ g^{-1}$ for some
$g:\sfB\to\Hom(\b{V},\b{V})$.

The aforementioned equivalence is given by $(\b{V},\nabla) \mapsto
\b{V}^{\nabla}$, where $\b{V}^{\nabla}$ is the local system, or
locally constant sheaf, of horizontal sections $\{\sigma \in
\Omega^0(\b{V})\mid \nabla(\sigma)=0\}$.  There is also a well
known equivalence between local systems on $\sfB$ and finite
dimensional complex representations of the fundamental group of
$\sfB$.  Note that isomorphic connections give rise to the same
representation.  Under these equivalences, the local systems
induced by the representations $\Phi^q(\b{t})$ and $\Psi^q(\b{t})$
correspond to flat connections on the vector bundles
$\b{K}^q\to\sfB$ and $\b{H}^q\to\sfB$, called Gauss-Manin
connections.

For weights $\bl$ that are both non-resonant and combinatorial, it
follows from \eqref{eq:NonResCohomology} that the only
non-vanishing cohomology vector bundle is $\b{H}^{\ell} \to \sfB$.
In \cite{T1}, Terao shows that this vector bundle is trivial, and
that the corresponding Gauss-Manin connection has logarithmic
poles along the irreducible components of the codimension one
divisor ${\sf D}=\overline\sfB \setminus \sfB$, where
$\overline\sfB$ denotes the closure of $\sfB$ in
$(\CP^{\ell})^{n}$.

For general weights, the vector bundles $\b{K}^q\to\sfB$ and
$\b{H}^q\to\sfB$ need not be trivial.  However, the restriction of
any one of these vector bundles to a circle is trivial, since any
map from the circle to the relevant classifying space is
null-homotopic.  Thus we make a local study of these bundles,
their Gauss-Manin
connections, and the corresponding local systems and fundamental
group representations as follows.

Let $\c\in\pi_{1}(\sfB,\sfb)$, and choose a representative path
$\tilde g:I\to\sfB$.  Then, of course, $\tilde g(0) = \tilde g(1) =
\sfb$, so $\tilde g$ defines a map $g:\bS^{1} \to \sfB$.  If $\phi$
denotes one of the representations $\Phi^q(\b{t})$ or $\Psi^q(\b{t})$,
there is an induced representation of $\pi_{1}(\bS^{1},1)=\langle
\zeta \rangle = \Z$, given by $\zeta \mapsto X$, where $X=\phi(\c)$.
Denote the matrix of $X$ by the same symbol.

Now let $\b{V}\to \sfB$ denote one of the vector bundles
$\b{K}^q\to\sfB$ or $\b{H}^q\to\sfB$, and let $g^{*}(\b{V}) \to
\bS^{1}$ be the induced vector bundle over the circle.  Pulling
back the relevant Gauss-Manin connection $\nabla$, we have a
corresponding connection $g^{*}(\nabla)$ on the bundle
$g^{*}(\b{V}) \to \bS^{1}$, which, as noted above, is necessarily
a trivial vector bundle. Specifying the flat connection
$g^{*}(\nabla)$ on this trivial bundle amounts to choosing a
logarithm, $Y$, of the matrix $X$ arising from the above
representation.  In summary:
\begin{prop}
\label{prop:log}
The connection matrix $Y$ satisfies $X = \exp(-2\pi\ii Y)$.
\end{prop}

The representations $\Phi^{\bul}(\b{t})$ and $\Psi^{\bul}(\b{t})$ are
induced by the universal representation $\Phi^\bul(\b{x}):
\pi_{1}(\sfB,\sfb) \to \Aut_{\L}^{}(K_{\L}^{\bul})$, see Theorems
\ref{thm:univend} and \ref{thm:InducedRep}.  We now define a
corresponding formal connection.  Recall from Theorem
\ref{thm:approx}.\ref{item:approx2} that the Aomoto complex
$(A^{\bul}_{R},\mu^{\bul})=(A^{\bul}_{R}(\sfb),\mu^{\bul}_{\sfb})$ is
chain equivalent to the linearization (at $\b{1}\in\T$) of the
universal complex $(K^\bul_\L,\D^\bul(\b{x}))$.  Choosing bases
appropriately, we can assume that the Aomoto complex is equal to this
linearization, see the proof of \cite[Thm.~2.13]{CO1}.

For each $\c\in\pi_1(\sfB,\sfb)$ and each $q$, let
$\Omega^q(\b{y})(\c)$ denote the linear term in the power series
expansion of $\Phi^q(\b{x})(\c)$ in $\b{y}$, where
$\b{x}=\exp(\b{y})$, that is $x_j=\exp(y_j)$ for $1\le j\le n$.
This defines a map (in fact, a representation)
\begin{equation} \label{eq:AomotoRep}
\Omega^q(\b{y}):\pi_1(\sfB,\sfb) \longrightarrow
\End_R(A^q_R).
\end{equation}
By construction, the entries of the matrix of $\Omega^q(\b{y})(\c)$
are linear forms in $y_1,\dots,y_n$, with integer coefficients, see
Theorem \ref{thm:approx}.\ref{item:approx1}.  We call the collection
$\Omega^\bul(\b{y})$ the {\em formal connection}.

\begin{thm} \label{thm:main1}
The eigenvalues of the formal connection are integral linear forms in
the variables $y_1,\dots,y_n$.  In other words, for each
$\c\in\pi_1(\sfB,\sfb)$ and each $q$, the eigenvalues of the formal
connection matrix $\Omega^q(\b{y})(\c)$ are integral linear forms in
$y_1,\dots,y_n$.
\end{thm}
\begin{proof}
Recall from Theorem \ref{thm:univaut} that the eigenvalues of
$\Phi^q(\b{x})(\c)$ are monomials of the form $x_1^{m_1}\cdots
x_n^{m_n}$.  Since $\Omega^q(\b{y})(\c)$ is the linear term in the
power series expansion of $\Phi^q(\exp(\b{y}))(\c)$ in $\b{y}$, the
result follows.
\end{proof}

\begin{thm} \label{thm:main1a}
%For any arrangement $\A$ and any system of weights $\bl$, if
Let $\bl$ be a system of weights, and let
$\b{t}=\exp(-2\pi\ii\bl)$.  For each $\c\in\pi_1(\sfB,\sfb)$, the
evaluation $\Omega^q(\bl)(\c)$ of the formal connection matrix
$\Omega^q(\b{y})(\c)$ at $\bl$ is a Gauss-Manin connection matrix
corresponding to the automorphism $\Phi^{q}(\b{t})(\c)$.
\end{thm}
\begin{proof}
Given $\c\in\pi_{1}(\sfB,\sfb)$, the endomorphism
$\Omega^{q}(\b{y})(\c)$ of $A^{q}_{R}$ is the linearization of the
automorphism $\Phi^{q}(\b{x})(\c)$ of $K^{q}_{\L}$.  Recall from
Theorem \ref{thm:TrivialAtOne} that $\Phi^{q}(\b{1})(\c)=\id$.  It
follows that $\Omega^{q}(\b{y})(\c)$ may be realized as a
logarithmic derivative of $\Phi^{q}(\b{x})(\c)$ at $\b{t}=\b{1}$.
This being the case, we have
\[
\Phi^{q}(\b{x})(\c) = \exp(\Omega^{q}(\b{y})(\c)),
\]
where $\b{x}=\exp(\b{y})$.  Since $\b{t}=\exp(-2\pi\ii\bl)$, the
specialization $\b{x}\mapsto\b{t}$ yields
\[
\Phi^{q}(\b{t})(\c) = \exp(\Omega^{q}(-2\pi\ii\bl)(\c)).
\]
Thus a Gauss-Manin connection matrix $Y(\c)$ satisfies $-2\pi\ii
Y(\c) = \Omega^{q}(-2\pi\ii\bl)(\c)$, see Proposition
\ref{prop:log}.  Now the entries of $\Omega^{q}(\b{y})(\c)$ are
linear forms in %the variables
$y_{1},\dots,y_{n}$.  Consequently, we have
$\Omega^{q}(-2\pi\ii\bl)(\c) = -2\pi\ii \Omega^{q}(\bl)(\c)$.
Therefore, the specialization $\b{y} \mapsto \bl$ yields the
Gauss-Manin connection matrix $Y(\c)=\Omega^{q}(\bl)(\c)$.
\end{proof}

These results leads to the following theorem, the main result of
this paper, which provides an affirmative answer to the question
of Terao stated in the Introduction.

\begin{thm} \label{thm:main}
For any arrangement $\A$ and any system of weights $\bl$, if
$\c\in\pi_1(\sfB,\sfb)$, then the eigenvalues of a corresponding
Gauss-Manin connection matrix in local system cohomology are
evaluations of linear forms with integer coefficients, and are
thus integral linear combinations of the weights.
\end{thm}
\begin{proof}
Let $\c\in\pi_1(\sfB,\sfb)$, and consider the corresponding
induced bundles over the circle as discussed above. By Theorem
\ref{thm:main1a}, the Gauss-Manin connection on the vector bundle
$\b{K}^q \to \bS^1$ is given by the matrix $\Omega^q(\bl)(\c)$ for
each $q$. By Theorem~\ref{thm:main1}, the eigenvalues of the
connection matrix $\Omega^{q}(\bl)(\c)$ are evaluations at $\bl$
of linear forms with integer coefficients for each $q$. Passage to
cohomology yields a connection matrix
$\overline\Omega^{q}(\bl)(\c)$, corresponding to the cohomology
representation $\Psi^{q}(\b{t})(\c)$, whose eigenvalues satisfy
the same condition.
\end{proof}

\begin{rem}
As noted above, for weights $\bl$ that are both non-resonant and
combinatorial, the only non-vanishing cohomology vector bundle
$\b{H}^\ell \to \sfB$, corresponding to the representation
$\Psi^\ell(\b{t})$, is trivial.  A trivialization is given by the
$\beta${\bf nbc} basis for the local system cohomology group
$H^\ell(\sfM_\sfb;\LL(\sfb))$, see \cite{FT,T1}.  In this context,
Terao \cite{T1} shows that the Gauss-Manin connection is determined by
a connection $1$-form $\sum d \log {\sf D}_j \otimes \nabla_j$, where
$\nabla_j \in \End_\C H^\ell(\sfM_\sfb;\LL(\sfb))$, each $d \log {\sf
D}_j$ denotes a $1$-form on $\overline\sfB$ with a simple logarithmic
pole along the irreducible component ${\sf D}_j$ of the divisor $\sf
D=\overline\sfB\setminus\sfB$, and the sum is over all such
irreducible components.  In cases where the codimension of
$\overline\sfB$ in $(\CP^{\ell})^{n}$ is small, the endomorphisms
$\nabla_{j}$ have been explicitly determined, by Aomoto-Kita \cite{AK}
in the codimension zero case, and by Terao \cite{T1} in the
codimension one case.  See \cite{OT2} for an exposition of these
results.

If $\c_j$ is a simple loop in $\sfB$ linking the component
${\sf D}_j$, then the endomorphisms $\nabla_j$ and
$\Omega_j=\overline\Omega^{\ell}(\bl)(\c_j)$ give rise to conjugate
automorphisms $\exp(-2\pi\ii \nabla_j)$ and
$\Psi^\ell(\b{t})(\c_j)=\exp(-2\pi\ii\Omega_j)$ of
$H^\ell(\sfM_\sfb;\LL(\sfb))$.  It follows that the connections on the
trivial vector bundle over the circle corresponding to $\nabla_j$ and
$\Omega_j$ are isomorphic.  By Theorem \ref{thm:main}, the eigenvalues
of the latter connection matrix are integral linear combinations of
the weights.
\end{rem}

\section{Combinatorial Connections}
\label{sec:GM2}

In this section, we investigate the combinatorial implications of the
formal connection defined on the Aomoto complex.  Recall that the
Aomoto complex is a universal complex,
$(A^{\bul}_{R}(\A),a_{\b{y}}\wedge)$, with the property that the
specialization $y_j \mapsto \la_j$ calculates the Orlik-Solomon
algebra cohomology $H^*(A^{\bul}(\A),a_\bl\wedge)$.

Let $\b{A}^q \to \sfB$ be the vector bundle over the moduli space
whose fiber at $\sfb$ is $A^q(\A_\sfb)$, the $q$-th graded component
of the Orlik-Solomon algebra of the arrangement $\A_\sfb$.  Given
weights $\bl$, the cohomology of the complex
$(A^{\bul}(\A_\sfb),a_\bl\wedge)$ gives rise to an additional vector
bundle $\b{H}^q(A) \to \sfB$ whose fiber at $\sfb$ is the $q$-th
cohomology group of the Orlik-Solomon algebra,
$H^q(A^{\bul}(\A_\sfb),a_\bl\wedge)$.  Like their topological
counterparts studied in the previous sections, these combinatorial
vector bundles admit {\em combinatorial connections}.

Fix a basepoint $\sfb \in \sfB$, and denote the Aomoto complex of
$\A_{\sfb}$ by simply $(A^{\bul}_{R},a_{\b{y}}\wedge)$.  As before,
let $\mu^{\bul}(\b{y})$ denote the boundary map with respect to a
given basis.  Recall from \eqref{eq:AomotoRep} that the formal
connection is comprised of maps
$\Omega^\bul(\b{y}):\pi_1(\sfB,\sfb) \to \End_R^{}(A^q_R)$.

\begin{prop} \label{prop:AomotoCommute}
For each $\c\in\pi_1(\sfB,\sfb)$, the endomorphisms
$\Omega^q(\b{y})(\c):A^{q}_{R}\to A^{q}_{R}$, $0\le q\le \ll$,
comprise a chain map $\Omega^\bul(\b{y})(\c)$ of the Aomoto complex
$(A^{\bul}_{R}(\sfb),\mu^\bul(\b{y}))$.
\end{prop}
\begin{proof}
For $\c\in\pi_1(\sfB,\sfb)$, we have an automorphism
$\Phi^\bul(\b{x})(\c)$ of the universal complex
$(K^\bul_\L,\D^\bul(\b{x}))$ by Theorems \ref{thm:univend} and
\ref{thm:univaut}.  Write $\Phi^q=\Phi^q(\b{x})(\c)$ and
$\D^q=\D^q(\b{x})$, and consider these maps as matrices with entries
in $\L$.  Then, for each $q$,
\[
\D^q \cdot \Phi^{q+1} = \Phi^q \cdot \D^q.
\]

Now make the substitution $\b{x}=\exp(\b{y})$, and denote power series
expansions in $\b{y}$ by $\D^q = \sum_{k\ge 0}\D^q_k$ and
$\Phi^q=\sum_{k\ge 0}\Phi^q_k$.  In this notation,
$\Omega^q(\b{y})(\c) = \Phi^q_1$.  Comparing terms of degree two in
the above equality, we obtain
\begin{equation} \label{eq:degree2}
\D^q_0 \cdot \Phi^{q+1}_2 + \D^q_1 \cdot \Phi^{q+1}_1 + \D^q_2 \cdot
\Phi^{q+1}_0 =\Phi^q_0 \cdot \D^q_2 + \Phi^q_1 \cdot \D^q_1 + \Phi^q_2
\cdot \D^q_0.
\end{equation}
By Remark \ref{rem:KdotAtOne}, we have $\D^q_0=0$.  By Theorem
\ref{thm:approx}.\ref{item:approx2}, the linearization of $\D^q$ is
equal to the boundary map of the Aomoto complex,
$\D^q_1=\mu^q(\b{y})$.  Also, Theorem~\ref{thm:TrivialAtOne} implies
that $\Phi^q_0=\id$ and $\Phi^{q+1}_0=\id$.  These facts, together
with \eqref{eq:degree2}, imply that $\mu^q(\b{y}) \cdot \Phi^{q+1}_1
=\Phi^q_1 \cdot \mu^q(\b{y})$.  In other words, $\Phi^\bul_1 =
\Omega^\bul(\b{y})(\c)$ is a chain map of the Aomoto complex.
\end{proof}

So we write $\Omega^\bul(\b{y}):\pi_1(\sfB,\sfb) \to
\End_R(A^\bul_R)$.  By Theorem \ref{thm:main1}, the eigenvalues of the
formal connection $\Omega^\bul(\b{y})$ on the Aomoto complex
$A^\bul_R$ are integral linear forms in $\b{y}$.  Using this fact and
the above Proposition, we obtain the following combinatorial analogue
of Theorem \ref{thm:main}.

\begin{thm} \label{thm:main2}
For any arrangement $\A$ and any system of weights $\bl$, the
eigenvalues of the combinatorial connection in Orlik-Solomon algebra
cohomology are evaluations of linear forms with integer coefficients,
and are thus integral linear combinations of the weights.
\end{thm}
\begin{proof}
For $\c \in \pi_1(\sfB,\sfb)$, the formal connection
$\Omega^\bul(\b{y})(\c)$ is a chain map on the Aomoto complex, which
induces upon specialization the Gauss-Manin connection in
Orlik-Solomon algebra cohomology.  The result follows.
\end{proof}

\section{An Example}
\label{sec:example}

We conclude by illustrating the results of the previous sections with
an explicit example.  Let $\A$ be the arrangement in $\C^2$ with
hyperplanes
\[
\begin{aligned}
H_1&=\{u_1+u_2=0\},\qquad &H_2&=\{2u_1+u_2=0\},\\
H_3&=\{3u_1+u_2=0\},      &H_4&=\{1+5u_1+u_2=0\}.
\end{aligned}
\]

\subsection{Universal Complexes}
We first record the universal complex $K^{\bul}_{\L}$ and the Aomoto
complex $A^{\bul}_{R}$ of $\A$.  The universal complex is equivalent
to the cochain complex of the maximal abelian cover of the complement
$M=M(\A)$.  For any $\b{t}\in(\C^*)^4$, the specializations at $\b{t}$
of the two complexes are quasi-isomorphic.  The latter complex may be
obtained by applying the Fox calculus to a presentation of the
fundamental group of the complement, see for instance \cite{CS1}.  A
presentation of this group is
\begin{equation} \label{eq:pres}
\pi_1(M) = \langle \c_1,\c_2,\c_3,\c_4 \mid
[\c_{3}\c_{1},\c_{2}],\ [\c_{1}\c_{2},\c_{3}],\ [\c_{i},\c_{4}]
\text{ for } i=1,2,3 \rangle,
\end{equation}
where $[\a,\beta]=\a\beta\a^{-1}\beta^{-1}$.  Using this presentation,
we obtain
\begin{equation*}
K^\bul_\L:\qquad \L \xrightarrow{\,\ \Delta^0\ } \L^4
\xrightarrow{\,\ \Delta^1\ } \L^5,
\end{equation*}
where, in matrix form, $\Delta^0 = \Delta^0(\b{x}) =
\left[\begin{matrix} x_1-1 & x_2-1 & x_3-1 & x_4-1\end{matrix}\right]$
and
\[
\Delta^1 = \Delta^1(\b{x}) = \left[\begin{matrix}
x_{3}-x_{2}x_{3} & 1-x_3          & 1-x_{4} & 0       & 0\\
x_1x_{3}-1       & x_{1}-x_{1}x_3 & 0       & 1-x_{4} & 0\\
1-x_{2}          & x_{1}x_2-1     & 0       & 0       & 1-x_{4}\\
0                & 0              & x_1-1   & x_2-1   & x_{3}-1
\end{matrix}\right].
\]

By Theorem \ref{thm:approx}.\ref{item:approx2}, the Aomoto complex
$A^{\bul}_{R}$ is the linearization of the complex $K^{\bul}_{\L}$.
Fixing the $\b{nbc}$-basis \cite[\S5.2]{OT2} for the Orlik-Solomon
algebra of $\A$ yields a corresponding basis for $A_{R}^{\bul}$.  With
respect to this basis, the Aomoto complex is given by
\begin{equation*}
A^\bul_R:\qquad R \xrightarrow{\,\ \mu^0\ } R^4
\xrightarrow{\,\ \mu^1\ } R^5,
\end{equation*}
where $\mu^0 = \mu^0(\b{y}) =
\left[\begin{matrix} y_1 & y_2 & y_3 & y_4\end{matrix}\right]$ and
\[
\mu^1 = \mu^1(\b{y}) = \left[\begin{matrix}
-y_{2}      & -y_3      & -y_{4} & 0      & 0\\
y_{1}+y_{3} & -y_3      & 0      & -y_{4} & 0\\
-y_{2}      & y_{1}+y_2 & 0      & 0      & -y_{4}\\
0           & 0         & y_1    & y_2    & y_{3}
\end{matrix}\right].
\]

\subsection{The Moduli Space and Related Bundles}
The moduli space of the arrangement $\A$ was studied in detail by
Terao \cite{T1}, see also \cite[Ex.~10.4.2]{OT2}.  This moduli space
may be described as
\begin{equation*}
\sfB =\sfB(\A) = \left\{
\begin{pmatrix}
z^1_0 & z^2_0 & z^3_0 & z^4_0 & 1 \\
z^1_1 & z^2_1 & z^3_1 & z^4_1 & 0 \\
z^1_2 & z^2_2 & z^3_2 & z^4_2 & 0
\end{pmatrix}
\Biggm|
\begin{aligned}
D_{i,j,k}&=0 \text{ if } \{i,j,k\}=\{1,2,3\}\\
D_{i,j,k}&\neq 0 \text{ if } \{i,j,k\}\neq\{1,2,3\}
\end{aligned}
\right\}.
\end{equation*}
Here, $(z^i_0 : z^i_1 : z^i_2)\in\CP^2$ for $1\le i \le 4$, and
$D_{i,j,k}$ denotes the determinant of the submatrix of the above
matrix with columns $i$, $j$, and $k$, for $1 \le i<j<k\le 5$.  This
moduli space is smooth, see \cite[Prop.~9.3.3]{OT2}.  Recall the fiber
bundle $\pi:\sfM \to \sfB$ of \cite[\S3]{T1}, with fiber
$\pi^{-1}(\sfb)=\sfM_{\sfb}=M(\A_{\sfb})$, the complement of the
arrangement $\A_{\sfb}$ combinatorially equivalent to $\A$.  The total
space of this bundle is given by
\[
\sfM = \{(\sfb,\b{u}) \in \sfB \times \C^{2} \mid \b{u} \in
\sfM_{\sfb}\}.
\]

For brevity, in \ref{subsec:UGM}--\ref{subsec:res} below, we calculate
various representation and connection matrices for a single element
$\a\in\pi_1(\sfB,\sfb_0)$, where $\sfb_0\in\sfB$ is the basepoint
(corresponding to $\A$) given below.  View $\bS^{1}$ as the set of
complex numbers of length one, and define $g:\bS^{1} \to \sfB(\A)$,
$s\mapsto g(s)$, by the following formula.
\[
\sfb_0=
\begin{pmatrix}
\,0 & 0 & 0 & 1 & 1 \\
\,1 & 2 & 3 & 5 & 0 \\
\,1 & 1 & 1 & 1 & 0
\end{pmatrix}
\qquad
g(s) =
\begin{pmatrix}
\,0 & 0 & 0 & 1 & 1 \\
\,\frac{3-s}{2} & \frac{3+s}{2} & 3 & 5 & 0 \\
\,1 & 1 & 1 & 1 & 0
\end{pmatrix}
\]
Note that $g$ is a loop based at $\sfb_0$ about the divisor defined by
$D_{1,2,5}=0$ in $\overline \sfB \setminus \sfB$, so represents an
element $\a$ of the fundamental group $\pi_{1}(\sfB,\sfb_0)$.  We will
determine the action of $\a$ on the universal complex $K_\L^\bul$.

For this, consider the induced bundle $g^{*}(\sfM)$, with total space
\[
E=\left\{\bigl(s,(\sfb,\b{u})\bigr) \in \bS^{1} \times (\sfB \times
\C^{2}) \mid g(s)=\sfb\text{ and } \b{u} \in \sfM_{\sfb}\right\},
\]
and projection $\pi'\bigl(s,(\sfb,\b{u})\bigr)=s$.  A similar bundle
over $\bS^{1}$ arises in the context of configuration spaces.  We
refer to \cite{Bi} as a general reference on configuration spaces and
braid groups.  Let $F_{n}(\C)=\{\b{v}\in \C^{n}\mid v_{i}\neq
v_{j}\text{ if }i\neq j\}$ be the configuration space of $n$ ordered
points in $\C$, the complement of the braid arrangement.  There is a
well known bundle $p:F_{n+1}(\C) \to F_{n}(\C)$, which admits a
section.  Writing $F_{n+1}(\C)=\bigl\{(\b{v},w)\in
F_{n}(\C)\times\C\mid w \in\C\setminus \{v_{j}\}\bigr\}$, the
projection is $p(\b{v},w)=\b{v}$.  The fiber of this bundle is
$p^{-1}(\b{v}) = \C\setminus\{v_{j}\}$, the complement of $n$ points
in $\C$.

Define $g_{1}:\bS^{1}\to F_{4}(\C)$ by
$g_{1}(s)=(\frac{3-s}{2},\frac{3+s}{2},3,4)$.  This loop represents
the standard generator $A_{1,2}$ of the pure braid group
$P_{4}=\pi_1(F_{4}(\C),\b{v}_0)$, the fundamental group of the
configuration space $F_4(\C)$, where $\b{v}_0=(1,2,3,4)$.  Let
$g_{1}^{*}(F_5(\C))$ be the pullback of the bundle $p:F_5(\C) \to
F_4(\C)$ along the map $g_1:\bS^{1}\to F_{4}(\C)$.  The bundle
$g_{1}^{*}(F_5(\C))$ has total space
\[
E_{1}=\left\{\bigl(s,(\b{v},w)\bigr) \in \bS^{1} \times (F_{4}(\C)
\times \C) \mid g_{1}(s)=\b{v}\text{ and } w \in
\C\setminus\{v_{j}\}\right\},
\]
and projection $p'\bigl(s,(\b{v},w)\bigr)=s$.  The two induced bundles
$g^*(\sfM)$ and $g_{1}^{*}(F_5(\C))$ are related as follows.  If
$\b{v}=g_{1}(s)$ and $w\in\C\setminus\{v_{j}\}$, it is readily checked
that the point $\b{u}=(-1,w)$ is in $\sfM_{g(s)}$, the fiber of
$g^{*}(\sfM)$ over $s\in\bS^{1}$.  This defines a map $h:E_{1}\to E$,
$\bigl(s,(\b{v},w)\bigr) \mapsto\bigl(s,(g(s),\b{u})\bigr)$, where
$\b{u}=(-1,w)$.  Checking that $\pi'\circ h = p'$, we see that
$h:g_{1}^{*}(F_5(\C)) \to g^*(\sfM)$ is a map of bundles.

\subsection{Universal Representations and Formal Connections}
\label{subsec:UGM} The fiber bundles $g^{*}_1(F_5(\C))$ and
$g^{*}(\sfM)$ admit compatible sections, induced by the section of
the configuration space bundle $p:F_5(\C)\to F_4(\C)$ and the
bundle map $h$ defined above. Consequently, upon passage to
fundamental groups, we obtain  the following commutative diagram
with split short exact rows.
\[
\begin{CD}
1 @>>> \pi_{1}(\C\setminus\{v_{j}\}) @>>> \pi_{1}(E_{1})
@>>> \pi_{1}(\bS^{1}) @>>>1\\
@.     @VV{h_{*}}V      @VV{h_{*}}V         @| \\
1 @>>> \pi_{1}(M) @>>> \pi_{1}(E)     @>>> \pi_{1}(\bS^{1}) @>>>1
\end{CD}
\]

Via the bundle map $h$, the fiber $\pi_{1}(\C\setminus\{v_{j}\})$ of
$g_{1}^{*}(\xi_{1})$ may be realized as the intersection of the line
$\{u_{1}=-1\}$ with the fiber of $g^{*}(\xi)$ in $\C^{2}$.  Thus, the
map $h_{*}:\pi_{1}(\C\setminus\{v_{j}\}) \to \pi_{1}(M)$ is the
natural projection of the free group on four generators,
$\pi_{1}(\C\setminus\{v_{j}\})=\F_{4}=\langle
\c_{1},\c_{2},\c_{3},\c_{4}\rangle$, onto the group $\pi_{1}(M)$ with
presentation \eqref{eq:pres}.  Let $\zeta$ denote the standard
generator of $\pi_1(\bS^1,1)$, mapping to $A_{1,2} \in
P_4=\pi_1(F_4(\C),\b{v}_0)$ and to $\a\in\pi_1(\sfB,\sfb_0)$ under the
homomorphisms induced by the maps $g_1$ and $g$.  The action of
$\zeta$ on the free group $\F_4$ coincides with that of $A_{1,2}$ on
$\F_4$, and is well known.  It is given by the Artin representation:
\[
\zeta(\c_{i})=A_{1,2}(\c_{i})=
\begin{cases}
\c_{1}^{}\c_{2}^{}\c_{i}^{}\c_{2}^{-1}\c_{1}^{-1}
&\text{if $i=1$ or $i=2$,}\\
\c_{i}
&\text{otherwise.}
\end{cases}
\]

By virtue of the commutativity of the above diagram, this action
descends to an action of $\a\in\pi_1(\sfB,\sfb_0)$ on $\pi_{1}(M)$
defined by the same formula.  The resulting action of $\a$ on the
universal complex $K_{\L}^{\bul}$---the universal representation---may
be determined using the Fox calculus, see for instance \cite{CS2} for
similar computations.  The action on $K^0_\L$ is trivial since $\a$
acts on $\pi_1(M)$ by conjugation.  The action on $K^{1}_{\L}$ is
familiar.  It is obtained by applying the Gassner representation to
the pure braid $A_{1,2}$.  We suppress the calculation of the action
of $\a$ on $K^{2}_{\L}$, and record only the result below.

Denote the universal representation and formal connection matrices
corresponding to $\a\in\pi_1(\sfB,\sfb_0)$ by $\Phi^q =
\Phi^q(\b{x})(\a)$ and $\Omega^q = \Omega^q(\b{y})(\a)$ respectively.
These matrices provide chain maps of the universal and Aomoto
complexes:
\begin{equation*}
\begin{CD}
\L @>{\,\ \Delta^0\ }>> \L^4 @>{\,\ \Delta^1\ }>> \L^5\\
@VV{\Phi^0}V            @VV{\Phi^1}V              @VV{\Phi^2}V\\
\L @>{\,\ \Delta^0\ }>> \L^4 @>{\,\ \Delta^1\ }>> \L^5
\end{CD}
\qquad \qquad \qquad
\begin{CD}
R @>{\,\ \mu^0\ }>> R^4 @>{\,\ \mu^1\ }>> R^5\\
@VV{\Omega^0}V      @VV{\Omega^1}V        @VV{\Omega^2}V\\
R @>{\,\ \mu^0\ }>> R^4 @>{\,\ \mu^1\ }>> R^5
\end{CD}
\end{equation*}
and are given by $\Phi^{0}=1$, $\Omega^{0}=0$,
\begin{alignat*}{3}
\Phi^1 &= {\hskip -2pt}\left[\begin{matrix}
1-x_1+x_1x_2 & 1-x_2 & 0 & 0 \\
x_1-x_1^2    & x_1   & 0 & 0 \\
0            & 0     & 1 & 0 \\
0            & 0     & 0 & 1
\end{matrix}\right] {\hskip -4pt},
&\Omega^1 &= {\hskip -2pt}\left[\begin{matrix}
y_2  & -y_2 & 0 & 0 \\
-y_1 & y_1  & 0 & 0 \\
0    & 0    & 0 & 0 \\
0    & 0    & 0 & 0
\end{matrix}\right]{\hskip -4pt},
\\
\Phi^2 &= {\hskip -2pt}\left[\begin{matrix}
x_1x_2 & 0 & 0            & 0     & 0 \\
x_2-1  & 1 & 0            & 0     & 0 \\
0      & 0 & 1-x_1+x_1x_2 & 1-x_2 & 0 \\
0      & 0 & x_1-x_1^2    & x_1   & 0 \\
0      & 0 & 0            & 0     & 1
\end{matrix}\right]{\hskip -4pt}, \
&\Omega^2 &= {\hskip -2pt}\left[\begin{matrix}
y_1+y_2 & 0 & 0    & 0    & 0 \\
y_2     & 0 & 0    & 0    & 0 \\
0       & 0 & y_2  & -y_2 & 0 \\
0       & 0 & -y_1 & y_1  & 0 \\
0       & 0 & 0    & 0    & 0
\end{matrix}\right]{\hskip -4pt}.
\end{alignat*}

\subsection{Non-Resonant Local Systems}\label{subsec:nonres}
Let $\bl=(\la_1,\la_2,\la_3,\la_4)$ be a system of weights in $\C^4$,
and $\b{t}=(t_1,t_2,t_3,t_4)$ the corresponding point in $(\C^*)^4$.
The induced local system $\LL$ on $M$ is non-resonant, and
$H^2(M;\LL)\simeq\C^{2}$, provided the rank of the matrix
$\Delta^1(\b{t})$ is equal to three.  If this is the case, then $\b{t}
\neq \b{1}$ and $\rank \Delta^0(\b{t})=1$.

Let $\Xi(\b{x}):\L^{5}\to\L^{2}$ and $\Upsilon(\b{y}):R^{5}\to R^{2}$
be the linear maps with matrices
\[
\Xi=\Xi(\b{x}) = \left[\begin{matrix}
x_4-1      & 0          \\
0          & x_4-1      \\
x_3-x_2x_3 & 1-x_3      \\
x_1x_3-1   & x_1-x_1x_3 \\
1-x_2      & x_1x_2-1
\end{matrix}\right]
\quad \text{and} \quad
\Upsilon=\Upsilon(\b{y}) = \left[\begin{matrix}
y_4     & 0    \\
0       & y_4  \\
-y_2    & -y_3 \\
y_1+y_3 & -y_3 \\
-y_2    & y_1+y_2
\end{matrix}\right].
\]
Note that $\Upsilon$ is the linearization of $\Xi$.  It is readily
checked that $\Xi \circ \D^{1} = 0$, $\Upsilon \circ \mu^{1} = 0$, and
that $\rank \Xi(\b{t}) = 2$ if $\b{t}\in (\C^{*})^{4}$ induces a
non-resonant local system $\LL$ on $M$.  Consequently, the projection
$\C^{5} \simeq K^{2} \twoheadrightarrow H^{2}(M;\LL) \simeq \C^{2}$
may be realized as the specialization at $\b{t}$ of the map $\Xi$.

Via $\Xi:K^{2}_{\L} \to \L^{2}$ and $\Upsilon:A^{2}_{R} \to R^{2}$,
the chain maps $\Phi^\bul:K^{\bul}_{\L} \to K^{\bul}_{\L}$ and
$\Omega^\bul:A^{\bul}_{R} \to A^{\bul}_{R}$ induce maps
$\overline\Phi(\b{x}):\L^2 \to \L^2$ and $\overline\Omega(\b{y}):R^2
\to R^2$, given by
\[
\overline\Phi = \overline\Phi(\b{x}) =
\left[\begin{matrix}
x_1x_2 & 0 \\
x_2-1  & 1
\end{matrix}\right]
\quad \text{and} \quad
\overline\Omega = \overline\Omega(\b{y}) =
\left[\begin{matrix}
y_1+y_2 & 0 \\
y_2     & 0
\end{matrix}\right].
\]
Specializing at $\b{t} \in (\C^{*})^{4}$ and $\bl \in \C^{4}$ yields
the representation matrix $\Psi^{2}(\b{t})(\a)=\overline\Phi(\b{t})$
and the corresponding Gauss-Manin connection matrix
$\overline\Omega^{2}(\bl)(\a)=\overline\Omega(\bl)$ in the cohomology
of the non-resonant local system $\LL$.  These matrices are
\[
\Psi^{2}(\b{t})(\a)=
\left[\begin{matrix}
t_1t_2 & 0 \\
t_2-1  & 1
\end{matrix}\right]
\quad \text{and} \quad
\overline\Omega^{2}(\bl)(\a)=
\left[\begin{matrix}
\la_1+\la_2 & 0 \\
\la_2       & 0
\end{matrix}\right].
\]
Up to a transpose, the latter recovers Terao's calculation of the
connection matrix corresponding to the divisor $D_{1,2,5}$, denoted by
$\Omega_{4}$ in \cite[Ex.~10.4.2]{OT2}.

\subsection{Resonant Local Systems}\label{subsec:res}
Now let $\LL$ be a non-trivial resonant local system on $M$.  Such a
local system corresponds to a point $\b{1} \neq \b{t} \in
(\C^{*})^{4}$ satisfying $t_{1}t_{2}t_{3}=1$ and $t_{4}=1$.  For each
such $\b{t}$, we have $H^{1}(M;\LL)\simeq\C$ and
$H^{2}(M;\LL)\simeq\C^{3}$.  Representation and Gauss-Manin connection
matrices corresponding to the loop $\a\in \pi_{1}(\sfB,\sfb)$ in
resonant local system cohomology may be obtained by methods analogous
to those used in the non-resonant case above.

Define $\Xi:\L^{5}\to\L^{3}$ and $\Upsilon:R^{5}\to R^{3}$ by
\[
\Xi = \left[\begin{matrix}
x_1x_{2}-1 & 0       & 0       \\
x_{2}-1    & 0       & 0       \\
0          & x_{2}-1 & 0       \\
0          & 1-x_1   & x_{3}-1 \\
0          & 0       & 1-x_{2}
\end{matrix}\right]
\quad \text{and} \quad
\Upsilon = \left[\begin{matrix}
y_1+y_{2} & 0    & 0    \\
y_{2}     & 0    & 0    \\
0         & y_2  & 0    \\
0         & -y_1 & y_3  \\
0         & 0    & -y_2
\end{matrix}\right].
\]
As before, $\Xi \circ \D^{1} = 0$, and $\Upsilon \circ \mu^{1} = 0$,
and $\Upsilon$ is the linearization of $\Xi$.  For each $\b{t}\in
(\C^{*})^{4}$ satisfying $t_{1}t_{2}t_{3}=1$ and $t_{4}=1$, we have
$\rank \Xi(\b{t}) = 3$.  So the projection $\C^{5} \simeq K^{2}
\twoheadrightarrow H^{2}(M;\LL) \simeq \C^{3}$ may be realized as the
specialization $\Xi(\b{t})$.

Via $\Xi:K^{2}_{\L} \to \L^{3}$ and $\Upsilon:A^{2}_{R} \to R^{3}$,
the chain maps $\Phi^\bul:K^{\bul}_{\L} \to K^{\bul}_{\L}$ and
$\Omega^\bul:A^{\bul}_{R} \to A^{\bul}_{R}$ induce $\overline\Phi:\L^3
\to \L^3$ and $\overline\Omega:R^3 \to R^3$, given by
\[
\overline\Phi = \left[\begin{matrix}
x_1x_2 & 0        & 0       \\
0      & x_{1}x_2 & 1-x_{3} \\
0      & 0        & 1
\end{matrix}\right]
\quad \text{and} \quad
\overline\Omega = \left[\begin{matrix}
y_1+y_2 & 0         & 0      \\
0       & y_{1}+y_2 & -y_{3} \\
0       & 0         & 0
\end{matrix}\right].
\]
Specializing yields the representation matrix $\Psi^{2}(\b{t})(\a)$
and the corresponding Gauss-Manin connection matrix
$\overline\Omega^{2}(\bl)(\a)$ in the second cohomology of the
resonant local system $\LL$.

Using the universal complex $K^\bul_\L$ and the conditions satisfied
by a point $\b{t} \in (\C^*)^4$ inducing the resonant local system
$\LL$, one can show that the representation matrix
$\Psi^{1}(\b{t})(\a)$ in first cohomology is given by
$\Psi^{1}(\b{t})(\a)=\left[t_1t_2\right]$.  The corresponding
Gauss-Manin connection matrix is, of course,
$\overline\Omega^{1}(\bl)(\a)=\left[\la_1+\la_2\right]$.

\begin{rem}
For this arrangement, every local system $\LL$ is combinatorial.
Given $\LL$, there are weights $\bl\in\C^4$ for which the local system
cohomology $H^*(M;\LL)$ is quasi-isomorphic to Orlik-Solomon algebra
cohomology $H^*(A^\bul,a_\bl\wedge)$.  Thus, for such weights, the
combinatorial connection matrices in the cohomology of the
Orlik-Solomon algebra coincide with the Gauss-Manin connection
matrices in local system cohomology computed above.
\end{rem}

\bibliographystyle{amsalpha}

\end{document}